\documentclass[a4paper,10pt]{amsart}
\usepackage[utf8]{inputenc}
\usepackage{amsthm}
\usepackage{amsmath}
\usepackage{bm}
\usepackage{amsfonts}
\usepackage{amssymb}
\usepackage{mathtools}
\usepackage{appendix}
\usepackage{mathrsfs}
\usepackage{setspace}
\usepackage{xcolor}
\usepackage{textgreek}
\usepackage{todonotes}
\usepackage[a4paper, left=3.5cm, right=3.5cm, top=4cm, bottom=4cm]{geometry}
\usepackage{thmtools} 

\mathtoolsset{showonlyrefs}

\usepackage[pdfdisplaydoctitle,colorlinks,breaklinks,urlcolor=blue,linkcolor=blue,citecolor=blue]{hyperref} 

\newcommand{\C}{\mathbb{C}}

\newcommand{\N}{\mathbb{N}}

\newcommand{\R}{\mathbb{R}}

\newcommand{\U}{\mathcal{U}}

\newcommand{\XX}{\mathcal{X}}
\newcommand{\x}{\mathbf{x}}

\DeclareMathOperator{\re}{Re}
\DeclareMathOperator{\im}{Im}

\DeclareMathOperator{\lip}{Lip}

\renewcommand{\epsilon}{\varepsilon}
\renewcommand{\setminus}{\smallsetminus}

\newcommand{\set}[1]{\left\{#1\right\}}
\newcommand{\pa}[1]{\left(#1\right)}
\newcommand{\bra}[1]{\left[#1\right]}

\newcommand{\norm}[1]{\left\|#1\right\|}

\newtheorem{thm}{Theorem}[section]
\newtheorem{definition}[thm]{Definition}

\newtheorem{lemma}[thm]{Lemma}
\newtheorem{proposition}[thm]{Proposition}

\newtheorem{hypothesis}{Hypothesis}

\theoremstyle{remark}
\newtheorem{remark}[thm]{Remark}

\numberwithin{equation}{section}

\newenvironment{acknowledgements}{%
  \begin{abstract}
}{%
  \end{abstract}
}

\title[Collapse and Burst of SQG Vortices]{Collapse and Burst of generalized Surface Quasi-Geostrophic point Vortices}
\author[F. Grotto]{Francesco Grotto}
  \address{Dipartimento di Matematica, Università di Pisa, Largo Bruno Pontecorvo 5, 56127 Pisa, Italia}
  \email{francesco.grotto(at)unipi.it}
\author[U. Pappalettera]{Umberto Pappalettera}
  \address{Fakult\"at f\"ur Mathematik, Universit\"at Bielefeld, D-33501 Bielefeld, Germany}
  \email{upappale(at)math.uni-bielefeld.de}
\date\today

\begin{document}

\begin{abstract}
We consider the generalized Surface Quasi-Geostrophic point vortices dynamics, and identify a sufficient condition implying existence of bursts out of (and collapses into) any given initial configuration of vortices.
The condition is related to the stability of the linearized dynamics around three vortices evolving in a self-similar fashion.
\end{abstract}

\maketitle

\section{Introduction}\label{sec:introduction}
The evolution of generalized Surface Quasi-Geostrophic (gSQG) point vortices on the plane $\R^2 \simeq \C$ is a Hamiltonian system of singular ODEs describing the motion of $N$ points $z_1,\dots,z_N\in \C$ with intensities $\xi_1,\dots,\xi_N\in\R\setminus\set{0}$, evolving according to
\begin{equation}\label{eq:freePV}
	\dot{\bar z}_j(t)= i c_\alpha \sum_{k\neq j}\xi_k \frac{|z_j(t)-z_k(t)|^{\alpha-2}}{z_j(t)-z_k(t)},
	\quad c_\alpha = - \left( 2^\alpha \Gamma(\alpha/2)^2 \sin(\alpha\pi/2) \right)^{-1},
\end{equation}
where the parameter $\alpha \in (0,3)$, $\alpha \neq 2$ dictates the interaction between vortices.

The system is defined in such a way that the empirical measure $\omega=\sum_{j = 1}^N\xi_j \delta_{z_j}$
is a weak solution to the two-dimensional gSQG equations 
\begin{equation*}
\begin{cases}
\partial_t \omega + u\cdot \nabla \omega =0,\\
u=\nabla^\perp (-\Delta)^{-\alpha/2} \omega.
\end{cases}
\end{equation*}
The special value $\alpha = 2$ corresponds to 2D Euler dynamics, whereas $\alpha=1$ corresponds to Surface Quasi-Geostrophic equations. Values of $\alpha \in (1,2)$ can be seen as interpolating between these two models.
We distinguish the local interaction ($\alpha<2$) and the nonlocal interaction ($\alpha>2$), the Euler system being in between. The case $\alpha>3$ corresponds to velocity fields increasing with the distance from the vortex, and is considered unphysical. Finally, the case $\alpha=3$ yields a velocity field independent of the distance from the vortex, and is a tricky threshold case. 
We refer to \cite{IwWa10,Badin2018,godard2021vortex,chen2024sufficient} for additional details.

Early investigations on Euler point vortices date back to the work of Helmholtz \cite{Helmholtz1858}, and have later found a mathematically rigorous treatment in a series of works by Marchioro and Pulvirenti \cite{MaPu83,Marchioro88,MaPu94}. Many subsequent investigations study the statistical mechanics of this system \cite{Ca92,Bo99,Ki12,GrRo20,GrLuRo24}. 
The case $\alpha \neq 2$ has gained popularity more recently: in addition to the aforementioned references, we mention \cite{godard2023holder,donati2023holder,donati2024existence}.

For the 2D Euler dynamics, the authors have identified in \cite{ARMAa} a sufficient condition (\autoref{hyp:A} below) such that the following holds true: For every initial configuration of $N-2$ distinct point vortices on the plane, $N\geq 3$, it is always possible to find $N$ vortices solution of \eqref{eq:freePV} on a small time interval $(0,T)$ that \emph{burst out} of the assigned initial configuration. By time reversal, up to changing the sign of the intensities $\xi_1,\dots,\xi_N$, one can furthermore exhibit a solution \emph{collapsing into} any designed configuration of $N-2$ vortices.
\autoref{hyp:A} is a condition linked to the eigenvalues of the linearized dynamics of \eqref{eq:freePV} around a self-similar solution with $N=3$, and roughly speaking corresponds to (or rather implies) some sort of ``stability" of self-similar solutions to \eqref{eq:freePV} under the influence of external perturbations. 
A similar kind of stability is also important to show patch-to-vortex approximation results \cite{zbarsky}.
It is also interesting to compare \cite{ARMAa} with the construction by Vishik \cite{vishik2018instability}, where non-uniqueness of solutions to the forced Euler equations with vorticity in $L^p$ is instead a consequence of instability of self-similar solutions.
One of the key contributions of \cite{ARMAa} is the verification of the existence of a three-vortices configuration satisfying \autoref{hyp:A}. 

In the present paper, we follow the same program and extend the results in \cite{ARMAa} to $\alpha \neq 2$.
Before moving on, we give a precise definition of what we mean by burst and collapse of configurations of point vortices.
For the sake of simplicity, we only present the case of threefold bursts and collapses, that is the one considered in this paper. 

\begin{definition} \label{def:burst}
Fix $N \in \N$, $N \geq 3$.
A configuration of $N-2$ pairwise distinct point vortices with positions $y_j$ and non-zero intensities $\zeta_j$,
\begin{equation*}
    y=(y_1,\dots, y_{N-2}) \in \C^{N-2}\setminus \triangle^{N-2},\quad
    (\zeta_1,\dots, \zeta_{N-2})\in \R\setminus\set{0}
\end{equation*}
admits a burst (respectively a collapse) of vortices if there exist $T>0$ and a solution $z:(0,T) \to \C^N$ (respectively $z:(-T,0) \to \C^N$) of \eqref{eq:freePV} such that $\zeta_j=\xi_j$ for $j=1,\dots,N-3$, $\zeta_{N-2} = \xi_{N-2} + \xi_{N-1} + \xi_{N}$, and:
\begin{equation*}
    \lim_{t\to 0} z_j(t)=
    \begin{cases}
        y_j &\mbox{if } j=1,\dots,N-3,\\
        y_{N-2} &\mbox{if } j=N-2,N-1,N.
    \end{cases}
\end{equation*}
\end{definition}

Proving the existence of bursts for gSQG point vortices assuming \autoref{hyp:A} follows, without much technical difficulties, from the same arguments of \cite{ARMAa}.
We briefly recall them in \autoref{sec:main}.
More precisely, we have the following:
\begin{thm} \label{thm:main}
    Fix $N \in \N$, $N \geq 3$. Let $\alpha\in (0,3)$ and assume \autoref{hyp:A} below.
    Then, every configuration $y \in \C^{N-2}\setminus \triangle^{N-2}$ admits a burst/collapse.
\end{thm}

On the other hand, exhibiting a configuration of three point vortices satisfying \autoref{hyp:A} is, in general, not trivial. 
It requires being able to prescribe the real part of the eigenvalues of a $4 \times 4$ matrix with entries depending in a complicated way upon the initial positions of the vortices, see the coefficients $L_{hk}$ in \autoref{lem:changeofcoor}.

We recall that \autoref{hyp:A} is satisfied in the Euler case $\alpha=2$, by the results of \cite{ARMAa}.
The ultimate goal of the present paper is to show the validity of \autoref{hyp:A} for values of $\alpha \neq 2$.
Here we prove that \autoref{hyp:A} holds true for SQG point vortices $\alpha =1$, see \autoref{thm:main} in \autoref{sec:conto}.
Moreover, we have numerical evidence that the same applies to $\alpha\in (1-\epsilon,2+\epsilon)$, for some $\epsilon>0$ sufficiently small. 
Outside of this interval, the sufficient condition fails. However, there is no reason not to believe that bursts and collapses do happen in this case, too. It would be interesting to find conditions weaker than \autoref{hyp:A} under which \autoref{thm:main} can be proved.

\section{Proof of \autoref{thm:main}} \label{sec:main}

The proof of \autoref{thm:main} follows closely the arguments of \cite{ARMAa}, so we opt to discuss only its key ideas making clear how they immediately extend to the more general case $\alpha\neq 2$, and refer to \cite{ARMAa} for full details. 

Let us recall that the dynamics \eqref{eq:freePV} preserves the \emph{center of vorticity}, \emph{energy} (the Hamiltonian function in conjugate coordinates $\re z_j, \xi_j \im z_j$) and the \emph{moment of inertia},
\begin{equation*}
    C=\sum_{j=1}^N \xi_j z_j,\quad 
    H=\frac{1}{c_\alpha}\sum_{j\neq k} \xi_j \xi_k |w_j-w_k|^{\alpha-2} = 0, \quad
	L=\sum_{j\neq k} \xi_j \xi_k |w_j-w_k|^2 = 0.
\end{equation*}
If $C=0$ we say that the evolution is \emph{centered}.
A \emph{relative equilibrium} to \eqref{eq:freePV} is a solution such that all distances $|z_j-z_k|$ are constant in time.
    
The proof of \autoref{thm:main} can be split into four main steps, which we discuss in the remainder of this section.
The time-inverting transformation
    \begin{equation*}
        t\mapsto -t,\quad \xi_i\mapsto -\xi_i,\, \quad i=1,\dots,N,
    \end{equation*}
    leaves the dynamics \eqref{eq:freePV} unchanged, so without loss of generality we can reduce ourselves to only discuss the existence of bursts of vortices.

\subsection{Step 1: Change of variables} \label{ssec:change_var}

We consider first the case of $N=3$ vortices, for which an explicit characterization of self-similar motion is available. The case $\alpha = 2 $ was already considered in \cite{ARMAa}, so we focus on the case $\alpha \neq 2$ hereafter. 


\begin{proposition} \label{prop:self}
    Let $\alpha \neq 2$ and let  $w:I \to \C^3\setminus \triangle^3$ be a maximal solution of \eqref{eq:freePV}, $N=3$, $I \subset \R$ an open time interval. Suppose $w$ is centered. Then the following are equivalent.
	\begin{itemize}
		\item[\emph{i})] The solution is \emph{self-similar},
		that is there exist distinct $a_1,a_2,a_3 \in \C$ and $Z:I \to \C \setminus\{0\}$ such that $w$ takes the form
		\begin{equation} \label{eq:aZ}
		w_j(t) = a_j Z(t), \quad j=1,2,3.
		\end{equation}
		\item[\emph{ii})] Either the solution is a relative equilibrium, or $H=L=0$.
		\item[\emph{iii})] There exist $a,b\in\R$ such that
		\begin{equation}\label{eq:asrelation}
		i c_\alpha \sum_{k \neq j}\xi_k\frac{|a_j-a_k|^{\alpha-2}}{a_j-a_k}= \bar{a}_j (a-i b),
		\qquad  j=1,2,3.
		\end{equation}
	\end{itemize}
	Moreover, when these conditions hold: 
\begin{itemize}
\item
if $a=0$, then the solution is a relative equilibrium and $I=\R$;
\item
if $a > 0$ (resp. $a<0$), then the solution is a burst (resp. collapse), and for some $t_0, \theta_0 \in\R$ it holds $I=(t_0,\infty)$ (resp. $I=(-\infty, t_0)$) and 
	\begin{equation}\label{eq:selfsimilar}
	Z(t)=((4-\alpha)a(t-t_0))^{\frac{1}{4-\alpha}} \,e^{i \left( \theta_0 + \frac{b}{(4-\alpha)a} \log (t-t_0) \right)}.
	\end{equation}
\end{itemize}
\end{proposition}

The latter follows from known computations, we refer the reader to \cite{Reinaud21} in the case $\alpha=1$ and we provide a sketch of the proof in the general case for completeness.

\begin{proof}
    Let $w_{jk}(t)=w_j(t)-w_k(t)$. 
    The equations of motion \eqref{eq:freePV} imply, by an explicit computation (see \cite[Eqs. (33),(41)]{Badin2018}), the equations of relative motion
    \begin{equation}\label{eq:relmotion}
        \frac{d}{dt} |w_{23}(t)|^2 =
        2(\alpha-2)c_\alpha A(t)\xi_1
        \pa{|w_{13}(t)|^{\alpha-4}-|w_{12}(t)|^{\alpha-4}},
    \end{equation}
    where $A(t)$ is the signed area of the triangle with vertices $w_1(t),w_2(t),w_3(t)$, 
    and $A(t)>0$ corresponds to counter-clockwise orientation of the vertices. 
    By cyclic permutation of indices, a similar expression holds for $|w_{12}|^2$ and $|w_{13}|^2$.

    Assuming item i) the conservation of
    \begin{equation*}
        H=c_\alpha^{-1}|Z(t)|^{\alpha-2}\sum_{j\neq k} \xi_j \xi_k |a_j-a_k|^{\alpha-2}, \quad
		L=|Z(t)|^2\sum_{j\neq k} \xi_j \xi_k |a_j-a_k|^2,
    \end{equation*}
    implies that either $|Z(t)|$ is constant or $H=L=0$, that is item ii). 
    Assume now item ii), and observe that a centered relative equilibrium must also be self-similar, so we can focus on the case $H=L=0$. The latter condition allows to
	express the quantity $\xi_1\xi_2|w_{12}|^2$ in two ways, namely
	\begin{align*}
	\xi_1\xi_2|w_{12}|^2 &= |w_{12}|^{4-\alpha} \xi_1\xi_2|w_{12}|^{\alpha-2} =
	|w_{12}|^{4-\alpha} \left(-\xi_1\xi_3|w_{13}|^{\alpha-2}-\xi_2\xi_3|w_{23}|^{\alpha-2} \right),
    \end{align*}
    and
    \begin{align*}
	\xi_1\xi_2|w_{12}|^2 &= -\xi_1\xi_3|w_{13}|^2 -\xi_2\xi_3|w_{23}|^2,
	\end{align*}
    from which, equating right-hand sides and rearranging,
	\begin{equation}\label{eq:bellefrazioni}
	\xi_1\frac{|w_{13}|^{\alpha-4}-|w_{12}|^{\alpha-4}}{|w_{23}|^2}
	=\xi_2\frac{|w_{12}|^{\alpha-4}-|w_{23}|^{\alpha-4}}{|w_{13}|^2}.
	\end{equation}
	
	By \eqref{eq:relmotion} and since we are excluding relative equilibrium,
	\begin{equation*}
	\frac{d}{dt}\pa{\frac{|w_{23}|^2}{|w_{13}|^2}}=0 \, \Leftrightarrow \,
	\xi_1 (|w_{13}|^{\alpha-4}-|w_{12}|^{\alpha-4})|w_{13}|^2
	=\xi_2 (|w_{12}|^{\alpha-4}-|w_{23}|^{\alpha-4})|w_{23}|^2,
	\end{equation*}
	where the condition on the right is exactly \eqref{eq:bellefrazioni}. Thus the ratio $|w_{23}|/|w_{13}|$ is constant in time.
	By cyclically permuting indices, we obtain the same result on the ratios $|w_{12}|/|w_{23}|$ and $|w_{13}|/|w_{12}|$, and the solution is self-similar.

    Assuming item iii) one can directly check that $w_j(t)=a_j Z(t)$, $j=1,2,3$, with $Z$ given by \eqref{eq:selfsimilar}, is a solution of \eqref{eq:freePV}. 
    Next we show that item i) implies item iii). Assume the former and write $Z(t)=r(t)e^{i \theta(t)}$, where $r(t) = |Z(t)|$ and $\theta(t)$ are real-valued. We set $0 \in I$ and $Z(0)=1$ without loss of generality, and write $a_{jk}=w_{jk}(0)=a_j-a_k$.
We notice that for every $j=1,2,3$ we have by i)
	\begin{align} \label{eq:thetadot}
	\frac{\dot{\bar{w}}_j(t)}{\bar{w}_j(t)} = \frac{\dot{r}(t)}{r(t)}-i\dot{\theta}(t) .
	\end{align}
	Since by assumption $w_j$ is a solution of \eqref{eq:freePV}, we also have
    \begin{align} \label{eq:tilde_a=a}
    \frac{\dot{\bar{w}}_j(t)}{\bar{w}_j(t)}
    &= 
    \frac{i c_\alpha}{{\bar{w}_j(t)}} \sum_{k\neq j}\xi_k \frac{|w_j(t)-w_k(t)|^{\alpha-2}}{w_j(t)-w_k(t)}
    \\ \nonumber
    &= 
    r(t)^{\alpha-4}\frac{i c_\alpha }{{\bar{a}_j}} \sum_{k\neq j}\xi_k \frac{|a_j-a_k|^{\alpha-2}}{a_j- a_k}.
    \end{align}
	By \eqref{eq:thetadot}, for every fixed $t \in I$ this last quantity must be a complex number independent of $j$, 
	that is exactly condition \eqref{eq:asrelation}.

For the second part of the statement, we argue as follows.
Taking the real part in equation \eqref{eq:tilde_a=a} and using \eqref{eq:thetadot}, we obtain
\begin{align*}
\frac{\dot{r}(t)}{r(t)} = r(t)^{\alpha-4} a 
\quad\Rightarrow\quad
\frac{d}{dt} r(t)^{4-\alpha} = (4-\alpha) a
\end{align*}
This immediately implies $r(t)=|Z(t)|$ constant if $a = 0$, and $w$ is a relative equilibrium. On the other hand, for $a \neq 0$, from the equation above we uniquely determine $$r(t)=((4-\alpha)a(t-t_0))^{\frac{1}{4-\alpha}}$$ with $t_0 = \frac{1}{a(\alpha-4)}$.
In order to determine the angle $\theta(t)$, we take the imaginary part of \eqref{eq:thetadot} and use \eqref{eq:asrelation} and the expression for $r(t)$ above to get
	\begin{equation*}
	\theta(t) = \theta_0 - \frac{b}{(4-\alpha)a} \log(t-t_0).\qedhere
	\end{equation*}
\end{proof}

We will study the stability of self-similar solutions to external perturbations by adopting the coordinates used in \cite{ARMAa}.
Fix distinct $a_1,a_2,a_3 \in \C$, $a_1 \neq 0$.
Introduce variables $r, \theta \in \R$ and $x_2, x_3 \in \C$ defined by
\begin{equation*}
(r e^{i \theta}, x_2, x_3)
=\pa{\frac{z_1}{a_1}, \frac{z_2}{z_1} - \frac{a_2}{a_1}, \frac{z_3}{z_1} - \frac{a_3}{a_1}},
\quad
z_1,z_2,z_3 \in \C.
\end{equation*}
The map $\Phi:(z_1,z_2,z_3) \mapsto (re^{i\theta},x_2,x_3)$ is a diffeomorphism of $(\C \setminus \{0\}) \times \C \times \C $ onto itself, and the diagonal set $\triangle^3$ is given by
\begin{equation*}
\triangle^3 = \set{x_2 = 1 - \frac{a_2}{a_1} }\cup \set{x_3 = 1 - \frac{a_3}{a_1} }\cup \set{x_2 - x_3 = \frac{a_3 - a_2}{a_1}}.
\end{equation*}
Let $w=w(t)$ be a self-similar solution to \eqref{eq:freePV} of the form \eqref{eq:aZ} with distinct $a_1,a_2,a_3 \in \C$, $a_1 \neq 0$, and such that $a > 0$ in \eqref{eq:asrelation} (that is, $w$ is a burst) and $t_0=0$ (without loss of generality, since the ODE system \eqref{eq:freePV} is autonomous).
In this coordinates $w(t)$ is given, for every $t \in I$, by
\begin{equation*}
x_2(t)=x_3(t)= 0,\quad 
r(t)=((4-\alpha)at)^{\frac{1}{4-\alpha}},
\quad 
\theta(t)=\theta_0+\frac{b}{(4-\alpha)a}\log t,
\end{equation*}
for some $\theta_0 \in \R$ and $a,b$ as in \eqref{eq:asrelation}. 

\subsection{Step 2: Including an external forcing term.} \label{ssec:external_forcing}

The coordinates $x_2=0$ and $x_3=0$ of the self-similar solution are equilibria, and \autoref{thm:main} is intimately linked to stability of this points with respect to external perturbations $f$ in the dynamics \eqref{eq:freePV}:
\begin{align} \label{eq:dynamics_f}
\dot{\bar z}_j(t)= i c_\alpha \sum_{k\neq j}\xi_k \frac{|z_j(t)-z_k(t)|^{\alpha-2}}{z_j(t)-z_k(t)}
+
f(t,z_j),
\quad
j=1,2,3,
\end{align} 
for some external field $f \in E_T = C([0,T),C^2(\C,\C)) \cap Lip([0,T),C(\C,\C))$. We identify $\C \simeq \R^2$, but differentiation will always be intended as real.
To study this stability, we rewrite the dynamics in \eqref{eq:dynamics_f} in terms of $r,\theta,x_2,x_3$.
In this new coordinates the dynamics is singular at $t=0$, but we can identify the leading order terms and effectively study the stability of the system around $x_2=x_3=0$ (which obviously is not an equilibrium for general $f \neq 0$).
In particular, we have the following:
\begin{lemma}\label{lem:changeofcoor}
There exist $C,\rho>0$ such that, for $r>0,\theta\in\R$ and $|x_2|, |x_3|, |x_2-x_3| < \rho$ it holds $(re^{i\theta},x_2,x_3) \notin \triangle^3$ and the system \eqref{eq:dynamics_f} is given in terms of $r,\theta,x_2,x_3$ by
	\begin{align*}
	\frac{d}{dt}(r^2)
	&=
	2ar^{\alpha-2} +\omega_r(x_2,x_3) +\frac{2}{|a_1^2|} \re \pa{z_1 f(t,z_1)},
	\\
	\frac{d}{dt}\theta
	&=
	\frac{b}{r^{4-\alpha}} + \frac{\omega_\theta(x_2,x_3)}{r^{4-\alpha}}+\im\pa{ \frac{\overline{f(t,z_1)}}{z_1} },
	\\
	\frac{d}{dt}x_j 
	&=
	\frac{L_{j}(x_2,x_3,\overline{x_2},\overline{x_3})}{r^{4-\alpha}} +
	\frac{\omega_j(x_2,x_3,x_2-x_3)}{r^{4-\alpha}}
	+\frac{1}{z_1} \overline{f(t,z_j)}
	-\frac{z_j}{z_1^2}\overline{f(t,z_1)},	
	\end{align*}
for $j=2,3$, where:
	\begin{itemize}
		\item $\omega_r,\omega_\theta:\set{z\in\C:|z|< \rho}^2\to\C$ are smooth functions such that 
		\begin{equation*}
		|\omega_{r,\theta}(x_2,x_3)| \leq C \left(|x_2| + |x_3|\right),\quad
		|\nabla \omega_{r,\theta}(x_2,x_3)| \leq C;
		\end{equation*}
		\item $\omega_2,\omega_3:\set{z\in\C:|z|< \rho}^3\to\C$ are smooth functions such that
		\begin{gather*}
		|\omega_{j}(x_2,x_3,x_2-x_3)| \leq C \left(|x_2|^2 + |x_3|^2\right),\\
		|\nabla \omega_{j}(x_2,x_3,x_2-x_3)|\leq C
		\left(|x_2| + |x_3|\right);
		\end{gather*} 
		\item $L_2,L_3:\C^4\to\C$ are $\C$-linear functions given by
		\begin{align*}
		L_2(x_2,x_3,\overline{x_2},\overline{x_3})
		&=
		(-a-ib)x_2 
		+
		L_{13}\overline{x_2}
		+
		L_{14}\overline{x_3},
		\\
		L_3(x_2,x_3,\overline{x_2},\overline{x_3})
		&= 
		(-a-ib)x_3
		+
		L_{23}\overline{x_2}
		+
		L_{24}\overline{x_3},
		\end{align*}
with $L_{13},L_{14},L_{23},L_{24} \in \C$ given by
\begin{align*}
L_{13} =&
-\frac{i c_\alpha}{|a_1|^2}
\overline{
	a_1^2 \xi_3 \left(
	(\alpha-2)|a_2-a_3|^{\alpha-4}-\frac{|a_2-a_3|^{\alpha-2}}{(a_2-a_3)^2} \right)}
\\
&\quad -\frac{i c_\alpha}{|a_1|^2}
\overline{
	a_1^2 \xi_1 \left(
	(\alpha-2)|a_2-a_1|^{\alpha-4}-\frac{|a_2-a_1|^{\alpha-2}}{(a_2-a_1)^2} \right)}
\\
&\quad -\frac{i c_\alpha}{|a_1|^2}
\frac{{a_2}}{{a_1}}
\overline{
	a_1^2 \xi_2 \left(
	(\alpha-2)|a_1-a_2|^{\alpha-4}-\frac{|a_1-a_2|^{\alpha-2}}{(a_1-a_2)^2} \right)
}, \\
L_{14} =&
-\frac{i c_\alpha}{|a_1|^2}
\overline{
	a_1^2 \xi_3 \left(
	(\alpha-2)|a_2-a_3|^{\alpha-4}-\frac{|a_2-a_3|^{\alpha-2}}{(a_2-a_3)^2} \right)} \\
&\quad-\frac{i c_\alpha}{|a_1|^2}
\frac{a_2}{a_1}
\overline{
	a_1^2 \xi_3 \left(
	(\alpha-2)|a_1-a_3|^{\alpha-4}-\frac{|a_1-a_3|^{\alpha-2}}{(a_1-a_3)^2} \right)},\\
L_{23} =&
-\frac{i c_\alpha}{|a_1|^2}
\overline{
	a_1^2 \xi_2 \left(
	(\alpha-2)|a_3-a_2|^{\alpha-4}-\frac{|a_3-a_2|^{\alpha-2}}{(a_3-a_2)^2} \right)} \\
&\quad-\frac{i c_\alpha}{|a_1|^2}
\frac{a_3}{a_1}
\overline{
	a_1^2 \xi_2 \left(
	(\alpha-2)|a_1-a_2|^{\alpha-4}-\frac{|a_1-a_2|^{\alpha-2}}{(a_1-a_2)^2} \right)},\\
L_{24} =&
-\frac{i c_\alpha}{|a_1|^2}
\overline{
	a_1^2 \xi_2 \left(
	(\alpha-2)|a_3-a_2|^{\alpha-4}-\frac{|a_3-a_2|^{\alpha-2}}{(a_3-a_2)^2} \right)}
\\
&\quad -\frac{i c_\alpha}{|a_1|^2}
\overline{
	a_1^2 \xi_1 \left(
	(\alpha-2)|a_3-a_1|^{\alpha-4}-\frac{|a_3-a_1|^{\alpha-2}}{(a_3-a_1)^2} \right)}
\\
&\quad -\frac{i c_\alpha}{|a_1|^2}
\frac{{a_3}}{{a_1}}
\overline{
	a_1^2 \xi_2 \left(
	(\alpha-2)|a_1-a_3|^{\alpha-4}-\frac{|a_1-a_3|^{\alpha-2}}{(a_1-a_3)^2} \right)
}.
\end{align*}
	\end{itemize}
\end{lemma}

The proof of the previous lemma only requires elementary (but cumbersome) calculations, and we omit it. It follows the lines of \cite[Lemma 2.7]{ARMAa}, to which we refer for details.
To close this paragraph, notice that by the chain rule
\begin{align*}
    \frac{d}{dt}(r^{4-\alpha})
	&=
    \frac{4-\alpha}{2} r^{2-\alpha} \frac{d}{dt}(r^2) 
    \\
    &=
	(4-\alpha)a  + \frac{4-\alpha}{2} r^{2-\alpha}\omega_r(x_2,x_3) + r^{2-\alpha} \frac{4-\alpha}{|a_1^2|} \re \pa{z_1 f(t,z_1)},
\end{align*}
in particular if $r=0$ at time $t=0$, then $r(t)^{4-\alpha} \sim (4-\alpha)at$ for small values of $t$ up to a small error.

\subsection{Step 3: Stability of self-similar solutions under external perturbations} \label{ssec:stability_external}
The key step in the construction of solutions leading to burst/collapse consists in finding a configuration of three free vortices $w_1,w_2,w_3$ on the plane, evolving in a self-similar fashion, that is sufficiently ``stable" to external perturbations, in a suitable sense which will become clear soon (cf. \autoref{hyp:A} below).
This notion of stability revolves around the change of coordinates $\Phi:(z_1,z_2,z_3) \mapsto (x_1,x_2,x_3)$ introduced in \autoref{ssec:change_var} above, and stability (around the equilibrium points $x_2=x_3=0$) is referred to a condition of the linearized dynamics under external perturbations given by \autoref{lem:changeofcoor}.

\begin{hypothesis}\label{hyp:A}
For every $\xi \neq 0$ there exists a choice of parameters $a_1,a_2,a_3 \in \C$, $\xi_1,\xi_2,\xi_3 \in \R$, such that \eqref{eq:asrelation} is satisfied with $a>0$, $\xi_1+\xi_2+\xi_3=\xi$ and the matrix
\begin{equation*}
L   = \pa{L_2,L_3,\bar L_2,\bar L_3}
=\pa{\begin{array}{cccc}
	-a-i b 	& 	0 	& L_{13} & L_{14}  	\\
	0 	& 	-a-i b 	& L_{23} & L_{24}  	\\
	\overline{L_{13}} & \overline{L_{14}} & -a + i b & 0						\\
	\overline{L_{23}} & \overline{L_{24}} & 0 	& -a + i b 
	\end{array}}
\end{equation*}
has eigenvalues whose real part is equal to $-a$.
\end{hypothesis}

Let us briefly comment on the relevance of the above condition on the eigenvalues in the arguments at the beginning of the section.
For notational convenience set
\begin{equation*}
    \x=(x_2,x_3,\overline{x_2},\overline{x_3}),\quad
    L(\x)=(L_2,L_3,\bar L_2,\bar L_3)(\x).
\end{equation*}
The equation for the time derivative of $\x$ reads
	\begin{align*}
	\frac{d}{dt}\x 
	&= 
	\frac{L(\x)}{(4-\alpha)at}  + \Xi,
	\end{align*}
where $\Xi = (\Xi_2,\Xi_3,\overline{\Xi}_2,\overline{\Xi}_3)$ is a remainder (in terms of magnitude at the singular time $t=0$) given by
	\begin{align*}
	\Xi_j &= L_j(\x)\pa{\frac1{r^{4-\alpha}}-\frac1{(4-\alpha)at}} +
	\frac{\omega_j(x_2,x_3,x_2-x_3)}{r^{4-\alpha}}
	+\frac{1}{z_1} \overline{f(t,z_j)}
	-\frac{z_j}{z_1^2}\overline{f(t,z_1)}.
	\end{align*}
Let us denote projections on single coordinates by $x_j=P_j[\x]$, $j=2,3$.
By Duhamel's formula, trajectories ${x}_j$ starting from ${x}_j=0$ at time $t=0$ can be represented as
\begin{equation*}
	{x}_j(t) = P_j \bra{\exp\left( \frac{\log t}{(4-\alpha)a} L \right)\int_{0}^{t}
		\exp\left( -\frac{\log s}{(4-\alpha)a} L \right)
		\Xi(s) ds}.
\end{equation*}
This makes it clear how the condition on the eigenvalues of $L$ allows good a priori estimates for $x_j$, 
making it possible to reinterpret the existence of a solution to the system of \autoref{lem:changeofcoor} as a fixed point problem. 
Indeed, if $L$ has all eigenvalues with real part equal to $-a$, then the matrix exponentials in the formula above satisfy the bounds (with respect to the operator norm)
\begin{align} \label{eq:expL}
\left\| \,\exp\left( \frac{\log t}{(4-\alpha)a} L \right)\right\|
\leq
t^{-\frac{1}{4-\alpha}},
\quad
\left\| \,\exp\left( -\frac{\log s}{(4-\alpha)a} L \right)\right\|
\leq
s^{\frac{1}{4-\alpha}}.
\end{align} 

In addition to this, \autoref{hyp:A} implies the existence of a self-similar solution $w$ (recall \autoref{prop:self} and the fact that $a_j,\xi_j$ satisfy \eqref{eq:asrelation}), that we fix hereafter. The condition of the eigenvalues further implies the following notion of stability for $w$.
Denote $\| f-g \|_{\infty} = \sup_{t \in (0,T)} \sup_{z\in \C} | f(t,z)-g(t,z)|$.

\begin{proposition} \label{prop:existenceexternal}
Assume \autoref{hyp:A}.
For all $M,\rho>0$ there exists $T^*>0$ such that, for any $T<T^*$ and $f \in E_T$ with $\|f\|_{E_T} \leq M$, there exists a solution $z:(0,T) \to \C^3 \setminus 
\triangle^3$ of \eqref{eq:dynamics_f} of class $C^1$ such that $\lim_{t \to 0} z_j(t)=0$.
In addition, there exists a constant $C$ such that 
\begin{align*}
    \sup_{t \in (0,T)} |z(t)| \leq \rho,
    \quad
    [z_j]_{C^\frac{1}{4-\alpha}_T} \leq C_\alpha[w_j]_{C^\frac{1}{4-\alpha}},
    \quad
    j=1,2,3.
\end{align*}
Moreover, if $z, z'$ are solutions respectively associated to the external forcing $f,g\in E_T$, with $\norm{f}_{E_T}, \norm{g}_{E_T} \leq M$, then it holds 
	\begin{align} \label{eq:fg}
	\sup_{\substack{0<t < T,\\j=1,2,3}} |z_j(t) - z_j'(t)| \leq C T^{{\frac{1}{4-\alpha}}} \|f-g\|_{\infty},
	\end{align}
	where $C>0$ is a constant depending only on $\xi,a,b,a_1,a_2,a_3,M, \alpha$.
\end{proposition}

The previous proposition implies the existence of a solution to the forced system \eqref{eq:dynamics_f} with $N=3$, bursting out of the origin at time $t=0$, continuous with respect to the external forcing $f$.
This corresponds to a notion of stability for $w$ (the particular self-similar solution whose existence is implied by \autoref{hyp:A}). 
To see this, take $g=0$ in \eqref{eq:fg}: a close inspection to the proof of the proposition gives $z'=w$, and so 
\begin{align*} 
	\sup_{\substack{0<t < T,\\j=1,2,3}} |z_j(t) - w_j(t)| \leq C T^{{\frac{1}{4-\alpha}}} \|f\|_{\infty}.
	\end{align*}
We point out that we do not claim that \emph{every} solution to \eqref{eq:dynamics_f} satisfies the same continuity with respect to the external forcing. 
Indeed, the solution constructed in \autoref{prop:existenceexternal} is obtained by fixing a particular ``angle'' (encoded in the condition $\eta_T = 0$ in the proof below).
In particular, if \eqref{eq:fg} held true for every solution to \eqref{eq:dynamics_f}, then one would have uniqueness
for the unforced system, which obviously does not hold since the choice of $\theta_0$ in \eqref{eq:selfsimilar} is arbitrary.
Our result should rather be seen as a continuous selection of solutions, and uniqueness is restored distinguishing how the burst out of the origin happens: that is, keeping track of the relative positions of vortices (encoded in the parameters $x_2,x_3$) and angle with respect to the horizontal axis (parameter $\theta$).

\begin{proof}[Proof of \autoref{prop:existenceexternal}]   
We first look for a solution of the system in \autoref{lem:changeofcoor} as a fixed point of a continuous map. 
Let $a_1,a_2,a_3$ and $a,b$ be as in \autoref{lem:changeofcoor}.
Besides those coefficients, the datum of the problem is the external field $f$, on which we can assume without loss of generality that $f(0,0)=0$, cf. \cite[Remark 2.9]{ARMAa}.

To set up the fixed point argument, take $T^*$ sufficiently small and for every $T < T^*$ define the convex, compact, non-empty space
\begin{align}\nonumber
\U_T &= \Big \{  (\zeta,\eta,x_2,x_3) \in C((0,T),\R_+ \times \R \times \C^2): \\
\label{eq:bc}
&\qquad \lim_{t \to T } \eta(t) = 0;\\
\label{eq:distance}
&\qquad |\zeta(t) - (4-\alpha)at| \leq t^{1+\frac{1}{4-\alpha}}, 
\,|x_2(t)|,|x_3(t)| \leq t \quad \forall t \in (0,T);\\
\label{eq:holder}
&\qquad \bra{\zeta}_{C^{\frac{1}{2}}_T},
\bra{\eta}_{C^{\frac{1}{2}}_T}, \bra{x_2}_{C^{\frac{1}{2}}_T}, \bra{x_3}_{C^{\frac{1}{2}}_T} \leq 1 \Big\},
\end{align}
and $\XX_T  = \Psi(\U_T)$, where the homeomorphism $\Psi: \U_T \to C((0,T),\C^3)$ is given by
\begin{gather*}
	\Psi(\zeta,\eta,x_2,x_3) = (r e^{i \theta},x_2,x_3),
	\\
 \quad r(t)^{4-\alpha} = \zeta(t), \quad \theta(t) = \eta(t) + \frac{b}{(4-\alpha)a} \log t.
	\end{gather*}

Let us assume $u=(\zeta,\eta,x_2,x_3)\in\mathcal{U}_{T}$  and rewrite equations of \autoref{lem:changeofcoor} in integrated form.
We have:
\begin{align}
	\nonumber
	\zeta(t) &= (4-\alpha)at + \int_0^t R(s,u_s) ds, \\
	\nonumber
	\eta(t)&= \int_{T}^t \Theta(s,u_s) ds,\\
	\label{eq:xequation}
	\x(t) &= \int_0^t \frac{L(\x_s)}{(4-\alpha)as} ds + \int_0^t \Xi(s,u_s) ds,
	\end{align}
	where 
	\begin{align*}
	R(s,u) 
    &=
    \frac{4-\alpha}{2}\zeta^{\frac{2-\alpha}{4-\alpha}}\omega_{r}(x_2,x_3) 
    +
    \frac{4-\alpha}{|a_1^2|} \zeta^{\frac{2-\alpha}{4-\alpha}}\re\pa{z_1 f(s,z_1)},
    \\
	\Theta(s,u) 
    &= 
    \left(\frac{b}{\zeta} -\frac{b}{(4-\alpha)as} \right)
	+ 
    \frac{\omega_{\theta}(x_2,x_3)}{\zeta}
	+
    \im\pa{ \frac{\overline{f(s,z_1)}}{z_1} },
    \\
	\Xi_j(s,u) 
    &= 
    L_j(\x)\pa{\frac1{\zeta}-\frac1{(4-\alpha)as}} +
	\frac{\omega_j(x_2,x_3,x_2-x_3)}{\zeta}
	+\frac{1}{z_1} \overline{f(s,z_j)}
	-\frac{z_j}{z_1^2}\overline{f(s,z_1)}.
\end{align*}

If $u \in \mathcal{U}_T$ is a solution of the system above, then $(r,\theta,x_2,x_3)=(\zeta^{\frac{1}{4-\alpha}},\eta+\frac{b}{(4-\alpha)a}\log t,x_2,x_3)$ solves the system of \autoref{lem:changeofcoor}.
Moreover, there exists a constant $C$ depending only on $a_i,\xi_i$ and $\|f\|_{E_T}$ such that for every solutions $u,u' \in \U_T$  
	\begin{align*}
	    |R(s,u_s)| &\leq C s,
     \\
     |R(s,u_s) - R(s,u'_s)| &\leq C \left(|\zeta_s - \zeta'_s|^{\frac{2-\alpha}{4-\alpha}}+|\zeta_s - \zeta'_s|^{\frac{1}{4-\alpha}} + |\eta_s - \eta'_s|+ |\x_s - \x'_s| \right),
	\\
	|\Theta(s,u_s)| &\leq C, 
    \\
    |\Theta(s,u_s) - \Theta(s,u'_s)| &\leq \frac{C}{s} \left(|\zeta_s -\zeta'_s|^{\frac{1}{4-\alpha}} + |\eta_s - \eta'_s|+ |\x_s - \x'_s| \right),
	\\
	|\Xi(s,u_s)| &\leq C s^{\frac{1}{4-\alpha}},
    \\
    |\Xi(s,u_s) - \Xi(s,u'_s)| &\leq C 
	\left(|\zeta_s - \zeta'_s|^{\frac{1}{4-\alpha}} + |\eta_s - \eta'_s|+ |\x_s - \x'_s| \right).
	\end{align*}	

In order to prove the existence of a solution $u=(\zeta,\eta,x_1,x_2)$ of the system above one can argue as in \cite[Lemma 2.12]{ARMAa}: there exists a $T^*$ sufficiently small such that, for any $T<T^*$, the relation $\Gamma(u)=\Gamma(\zeta,\eta, x_2, x_3) = (\tilde{\zeta},\tilde{\eta}, \tilde{x}_2, \tilde{x}_3)$ given by
	\begin{align*}
	\tilde{\zeta}(t) &= (4-\alpha)at + \int_0^t R(s,u_s) ds, \\
	\tilde{\eta}(t) &= \int_{T}^t \Theta(s,u_s) ds,\\
	\tilde{x}_j(t) &= P_j \bra{\exp\left( \frac{\log t}{(4-\alpha)a} L \right)\int_{t_0}^{t}
		\exp\left( -\frac{\log s}{(4-\alpha)a} L \right)
		\Xi(s,u_s) ds},
	\end{align*}
	defines a continuous map $\Gamma:\U_{T} \to \U_{T}$.
We point out that stability of the self-similar solution plays a crucial role in the proof of continuity of $\Gamma$, since the bounds in \eqref{eq:expL} imply
\begin{align*}
  |\tilde{x}_j(t)-\tilde{x}'_j(t)|
  &\leq
  t^{-\frac{1}{4-\alpha}} \int_0^t s^{\frac{1}{4-\alpha}} |\Xi(s,u_s)-\Xi(s,u'_s)| ds
  \\
  &\leq
  Ct \left( \|u-u'\|_{\infty} + \|u-u'\|_{\infty}^{\frac{1}{4-\alpha}} \right).
\end{align*}
 
 Moreover, up to choosing $T^*$ sufficiently small, any $r=\tilde{\zeta}^{\frac{1}{4-\alpha}}$ and $\theta=\tilde{\eta}+\frac{b}{(4-\alpha)a}\log t$, where $(\tilde{\zeta},\tilde{\eta},\tilde{x_1},\tilde{x_2})=\Gamma(u)$ with $u \in \mathcal{U}_T$, satisfies $re^{i\theta} \in C^{\frac{1}{4-\alpha}}((0,T),\C)$ and $[re^{i\theta}]_{C^{\frac{1}{4-\alpha}}_T} \leq C [Z]_{C^{\frac{1}{4-\alpha}}}$.
 Therefore, the first part of \autoref{prop:existenceexternal} descends from Schauder fixed point Theorem and change of coordinates.
Finally, continuous dependence on the external field follows the lines of \cite[Proposition 2.15]{ARMAa}.

\end{proof}

\begin{remark}
    In view of previous inequalities, more relaxed conditions on the eigenvalues $\lambda$ of the matrix $L$ (as for example $-C_1 a \leq Re(\lambda) \leq -C_2 a$ for suitable constants $C_1$ and $C_2$) should still be sufficient to prove \autoref{thm:main}. However, we have not found such weaker condition easier to verify in practice.
\end{remark}

\subsection{Step 4: Implications of stability}
Suppose we are given a suitable stable self-similar solution $w_j$, $j=1,2,3$, and $N-2$ points on the plane $y_1,\dots,y_{N-2}$. 
Placing the vortices $w_j$ such that their center of vorticity equals $y_{N-2}$, and assuming that vortices $y_1,\dots,y_{N-3}$ are well separated from the initial three for all times (this is always possible up to choosing a time $T$ sufficiently small), a fixed point argument allows to consider the influence of $y_1,\dots,y_{N-3}$ on $w_j$ as an external perturbations, and existence of a burst of vortices follows from the peculiar stability properties mentioned above. 
 
Let us introduce some notation in order to distinguish the bursting vortices from the others. 
By translation invariance, one can suppose without loss of generality $y_{N-2} = 0$.
Let $y(0) \in \C^{N-3}\setminus \triangle^{N-3}$ be a configuration of $N-3$ vortices with 
$y_k(0) \neq 0$ for every $k=1,\dots,N-3$ and intensities $\zeta_k \in \R$. 
In the origin $0 \in \C$ we consider a point vortex with intensity 
$\xi \in \R$, $\xi \neq 0$, which is the one we are going to split into three vortices 
of intensities $\xi_1,\xi_2,\xi_3$ as in \autoref{hyp:A}. 
In particular, \autoref{thm:main} descends from the following:

\begin{proposition} \label{thm:Nvortices}
	In the same setting as above, there exists $T>0$ and a solution $(z,y):(0,T) \to \C^3 \times \C^{N-3}$ of the full system:
	\begin{align} \label{eq:Nvorticesz}
	\qquad\dot{\bar{z}}_j	
	&= i c_\alpha  \sum_{\ell \neq j} \xi_\ell\frac{|z_j-z_\ell|^{\alpha-2}}{z_j - z_\ell}
	+i c_\alpha \sum_{h = 1}^N \zeta_h\frac{|z_j-y_h|^{\alpha-2}}{z_j - y_h}, \quad j=1,2,3,	\\
	\label{eq:Nvorticesy}
	\dot{\bar{y}}_k	
	&= i c_\alpha \sum_{\ell=1}^3 \xi_\ell\frac{|y_k - z_\ell|^{\alpha-2}}{y_k - z_\ell}
	+ i c_\alpha \sum_{h \neq k} \zeta_h\frac{|y_k - y_h|^{\alpha-2}}{y_k - y_h},
	\quad k=1,\dots,{N-3}, 
	\end{align}
	such that $\lim_{t\to 0} z_j(t) = 0$, $j=1,2,3$, and $\lim_{t \to 0} y_k(t) =y_k(0)$, $k=1,\dots,{N-3}$.
\end{proposition}

To set up the fixed point argument, recall the definitions of $\U_T$, $\XX_T= \Psi(\U_T)$ from \autoref{ssec:stability_external}.
Moreover, for $\rho>0$ sufficiently small introduce
\begin{align}\nonumber
\mathcal{Y}_{\rho,y_0,T} &= \left\{ y \in \lip((0,T),\C^{N-3}):\right.\\
\label{eq:bc_y}
&\qquad \lim_{t \to 0 } y_k(t) = y_k(0);\\
\label{eq:distance_y}
&\qquad |y_k(t) - y_k(0)| \leq \rho \quad \forall t \in (0,T);\\ 
\label{eq:distance_y0}
&\qquad |y_k(t)| \geq \rho \quad \forall t \in (0,T);\\ 
\label{eq:holder_y}
&\qquad \left.
\bra{y}_{\lip} \leq
|c_\alpha| (N-1) \rho^{\alpha-3} \max_{j,\ell} \{|\xi_j|,|\zeta_\ell|\} \quad
\right\}.
\end{align}

The sets $\U_T$ and $\mathcal{Y}_{\rho,y_0,T}$ are convex and compact with respect to the topology induced by the supremum norm, and non-empty since $\Psi^{-1}(w) \in \U_T$ and $y_k(t)\equiv y_k(0) \in \mathcal{Y}_{\rho,y_0,T}$. 

In order to prove \autoref{thm:Nvortices}, we regard our equations
as a fixed point of a continuous function on the space $\XX_T \times \mathcal{Y}_{\rho,y_0,T}$, and apply Schauder fixed point Theorem. 
We denote by
\begin{align} \label{eq:fNvortices}
f(t,p) = F(y_1(t),\dots,y_{N-3}(t);p) =
ic_\alpha \sum_{k=1}^{N-3}
\zeta_k\frac{|p-y_k(t)|^{\alpha-2}}{p-y_k(t)},
\quad p \in \C,
\end{align}
the field acting on $z_j$ due to the presence of vortices $y_k$.
As already done in \cite{ARMAa}, one can modify $F$ smoothly outside of the ball $B_\rho^{N-3} \times B_\rho$, $\rho>0$ sufficiently small, in such a way that for every $y \in \mathcal{Y}_{\rho,y_0,T}$ it holds $f \in E_T$ and $\norm{f}_{E_T} \leq M$, 
with $M$ depending only on $N$, $\xi_j$, $j=1,2,3$, $\zeta_k$, $k=1,\dots,N-3$ and $\rho$, without affecting the dynamics of vortices $z_j$ for $T<T^*$.

Analogously, we introduce the vector field acting on $y_k$, $k=1,\dots,{N-3}$:
\begin{align*}
g(t,p) = G(z_1(t),z_2(t),z_3(t);p) =
ic_\alpha \sum_{j=1}^3
\xi_j\frac{|p-z_j(t)|^{\alpha-2}}{p-z_j(t)},
\quad p \in \C.
\end{align*}

Next, define the map $\Gamma^\mathcal{X}: \mathcal{Y}_{\rho,y_0,T} \to \XX_T$ 
that to every $y \in \mathcal{Y}_{\rho,y_0,T}$ associates the unique solution $x \in \XX_T$ (after identification by change of coordinates) of \eqref{eq:dynamics_f} given by \autoref{prop:existenceexternal}, 
with external field $f = F(y)$ given by \eqref{eq:fNvortices} above.
Let now $\Gamma^\mathcal{Y}: \XX_T \times \mathcal{Y}_{\rho,y_0,T} \to C((0,T),\C^{N-3})$ 
be the map that associates to every $(x,y) \in \XX_T \times \mathcal{Y}_{\rho,y_0,T}$
the curve $\tilde{y}$ of components
\begin{align*}
\tilde{y}_k(t) =
\int_0^t \left( ic_\alpha\sum_{\ell \neq k}
\zeta_\ell\frac{|y_k-y_\ell|^{\alpha-2}}{y_k-y_\ell}  + G(\Phi^{-1}(x(s));y_k) \right) ds.
\end{align*}

Then, a solution to \eqref{eq:Nvorticesz} and \eqref{eq:Nvorticesy} is produced as a fixed point of $\Gamma(x,y) = (\Gamma^\mathcal{X}(y),\Gamma^\mathcal{Y}(x,y))$. 
The key property allowing us the apply Schauder Theorem is  the following:
    For $\rho,T>0$ sufficiently small, the map $\Gamma^\mathcal{Y}$ takes values in $\mathcal{Y}_{\rho,y_0,T}$ and the map $\Gamma:\XX_T \times \mathcal{Y}_{\rho,y_0,T} \to \XX_T \times \mathcal{Y}_{\rho,y_0,T}$ is continuous with respect to the supremum norm.

Given the preparations above, we can finally give the:
\begin{proof}[Proof of \autoref{thm:Nvortices}]
	Take $\rho$ and $T$ as above.
	Since $\Gamma:\XX_T \times \mathcal{Y}_{\rho,y_0,T} \to \XX_T \times \mathcal{Y}_{\rho,y_0,T}$ is continuous and $\XX_T \times \mathcal{Y}_{\rho,y_0,T}$ is non-empty and homeomorphic to a compact, convex subset of the Banach space $C((0,T),\C^3 \times \C^{N-3})$, Schauder fixed point Theorem yields existence of a fixed point $(x,y) \in \XX_T \times \mathcal{Y}_{\rho,y_0,T}$ for $\Gamma$. Let $z = \Phi^{-1}(x)$. Then $(z,y)$ is a solution of \eqref{eq:Nvorticesz}, \eqref{eq:Nvorticesy} satisfying $\lim_{t\to 0} z_j(t) = 0$, $j=1,2,3$, and $\lim_{t \to 0} y_k(t) =y_k(0)$, $k=1,\dots,N-3$, concluding the proof.
\end{proof}

\section{Parameters Satisfying \autoref{hyp:A}} \label{sec:conto}

We shall now discuss the existence of configurations for which \autoref{hyp:A} is satisfied. The problem can be reduced to checking algebraic relations. Given a value of $\alpha\in (0,3)$ one can exhibit a specific configuration satisfying \autoref{hyp:A}, and we will prove for instance:
\begin{thm} \label{thm:main}
    If $\alpha=1$ there exists a configuration of three point vortices with positions $a_1,a_2,a_3\in\C$ and intensities $\xi_1,\xi_2,\xi_3\in\R^\ast$ such that \autoref{hyp:A} is satisfied.
\end{thm}
A specific configuration can not satisfy \autoref{hyp:A} for all $\alpha\in (0,3)$, and it is not clear whether it is possible to provide an analytic expression for the involved parameters as functions of $\alpha$ for which \autoref{hyp:A} holds true. Therefore, we derive sufficient conditions for \autoref{hyp:A} in the case $\alpha=1$ and provide compelling numerical evidence that they can be satisfied for $\alpha$ in an open interval containing $[1,2]$.

By invariance of \eqref{eq:freePV} under orientation-preserving isometries of $\R^2$, in the remainder of this section we assume without loss of generality that $a_1 = 1/2$ and $a_2 = -1/2$.
We further restrict ourselves to the case $\xi_1=\xi_2=1$, so that the remaining free parameters are $(\xi_3,a_3)\in \R\times\C$. It is convenient to formulate the problem in terms of the side lengths
\begin{equation*}
    x=|a_{13}|, \quad y=|a_{23}|.
\end{equation*}

\autoref{hyp:A} requires the condition \eqref{eq:asrelation},
characterizing configurations giving rise to self-similar evolution, with $a>0$, that is excluding relative equilibria.
The condition $H=L=0$ is equivalent to \eqref{eq:asrelation} by \autoref{prop:self}, and with the above choices it becomes
\begin{equation}\label{eq:paramreduce}
		\xi_3= -\frac{1}{x^{\alpha-2} + y^{\alpha-2}}=- \frac{1}{x^{2} + y^{2}}>-2.
\end{equation}
The last inequality, which is required for the total intensity to be positive, follows from triangular inequality $x+y>|a_{12}|=1$ (unless positions $a_1,a_2,a_3$ are collinear, a relative equilibrium that we exclude) and AM-QM inequality.
Notice that \eqref{eq:paramreduce} has no positive solutions $(x,y)$ with $x,y$ both larger than $1$: By symmetry we can restrict to $x\in (0,1]$.
We can also write
\begin{equation} \label{eq:paramreduce2}
    a_3= \frac{x^2-y^2}{2} \pm i \sqrt{y^2-\frac{(x^2-y^2-1)^2}{4}}.
\end{equation}
In order to exclude collinear configurations, we impose $\im(a_3)\neq 0$ (that is, we ask that the argument of the square root be positive).
 This restricts the range of parameters $x,y$.
Overall, the first requirement of \autoref{hyp:A} translates to
\begin{equation}\label{eq:intermsofxy}
    \begin{cases}
        0<x<1,\, y>1-x,\\
        x^{\alpha-2}+y^{\alpha-2}=x^2+y^2,\\
        4y^2> (x^2-y^2-1)^2,
    \end{cases}
\end{equation}
and we are left to find choices of $x,y$ satisfying the above conditions and for which the stability part of \autoref{hyp:A} also holds true.

The choice of the sign of $\im(a_3)$ determines whether we deal with a burst or collapse. 
In order to have a burst, corresponding to $a>0$ as required by \autoref{hyp:A}, we shall choose $a_3$ with $\im(a_3)<0$.

Imposing that the eigenvalues of the matrix $L=L(x,y)$ have the form $-a + i \mu$, $\mu \in \R$, leads to an equation for $\mu$,
\begin{align*}
\mu^4 + \mu^2 (c_1 - 2b^2) + b^4 - b^2 c_1 + c_2  = 0,
\end{align*}
whose coefficients 
\begin{gather*}
    c_1 =L_{23}\overline{L_{14}} + 
L_{24}\overline{L_{24}} +
L_{13}\overline{L_{13}} +
L_{14}\overline{L_{23}}, \\
c_2 =L_{13}\overline{L_{13}}L_{24}\overline{L_{24}}+
L_{23}\overline{L_{23}}L_{14}\overline{L_{14}} -
L_{14}\overline{L_{13}}L_{23}\overline{L_{24}} -
L_{13}\overline{L_{14}}L_{24}\overline{L_{23}},
\end{gather*}
can be expressed in terms of $x,y$ by \autoref{lem:changeofcoor}, and $a,b$ can be expressed in terms of $x,y$ by \eqref{eq:asrelation}.
If the above equation for $\mu$ has 4 distinct real roots, then $L$ must have 4 distinct eigenvalues of the form $\lambda_j = -a+i\mu_j$, $j=1,\dots,4$, thus satisfying the second requirement of \autoref{hyp:A}. 
We therefore add to \eqref{eq:intermsofxy} the requirement 
\begin{equation}\label{eq:conditionmu}
2b^2 - c_1 \pm \sqrt{(2b^2 - c_1)^2 - 4(b^4 - b^2 c_1 + c_2)} > 0.
\end{equation}

\begin{proof}[Proof of \autoref{thm:main}]
    For $\alpha=1$, the second line of \eqref{eq:intermsofxy} is a cubic polynomial in $y$ with coefficients depending on $x$, namely
    \begin{align} \label{eq:cubic}
    y^3+\frac{x^3-1}{x}y-1
    =: y^3+p(x) y + q(x)
    =0.
    \end{align}
Next, observe that for every $x > 1/2$ it holds
    \begin{align*}
    \frac{q^2}{4} + \frac{p^3}{27} 
    =
    \frac{1}{4} - \frac{(1-x^3)^3}{27x^3}
    > 0.05150 > 0.
    \end{align*}
    As such, the solutions $y=y(x)$ of \eqref{eq:cubic} are explicitly given by Cardano's Formula, and the only real solution is 
    \begin{equation*}
    y=
    \sqrt[3]{{\frac12}+\sqrt{\frac{1}{4} - \frac{(1-x^3)^3}{27x^3}}}
    +
    \sqrt[3]{{\frac12}-\sqrt{\frac{1}{4} - \frac{(1-x^3)^3}{27x^3}}}.
    \end{equation*}
    The (exact) value $x=0.70190$ gives $y\simeq 1.30353$, and it is immediate to verify that this choice of $(x,y)$ fulfills all the conditions listed in \eqref{eq:intermsofxy}. We fix this choice of $(x,y)$ hereafter.
   By \eqref{eq:paramreduce} and \eqref{eq:paramreduce2}, paying attention to picking $a_3$ with $\im(a_3)<0$, we obtain the following approximate values for the vorticity and position of the third vortex:  
   \begin{equation*}
    \xi_3 \simeq -0.45623,
    \quad
    a_3 \simeq -0.60326 - 0.69426 i.
    \end{equation*}   
    Direct computations reveal that the configuration of vortices $a_1,a_2,a_3$ above, with respective intensities $\xi_1,\xi_2,\xi_3$, satisfies both \eqref{eq:intermsofxy} and \eqref{eq:conditionmu}, with
    \begin{equation*}
    a \simeq \frac{1}{\pi} \im  \left( 1 - 0.45623 \times \frac{1.42470}{1.10326 + 0.69426 i} \right)  \simeq  0.08453 > 0.\qedhere
    \end{equation*}
\end{proof}

\begin{figure}[t]
    \centering
        \label{fig:enter-label}
    \centering
    \includegraphics[width=10cm,height=9.42cm]{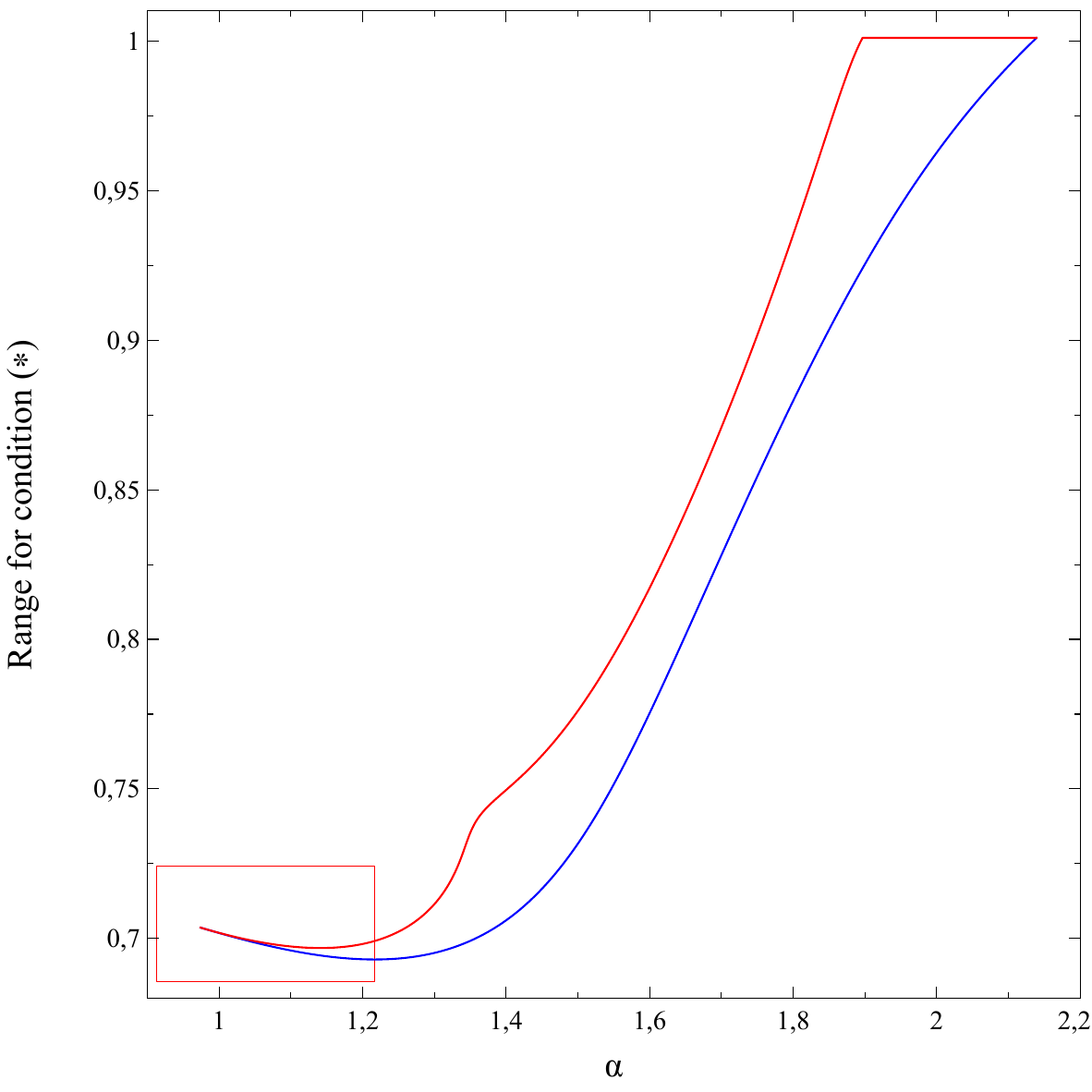}
    \centering
    \includegraphics[width=10cm,height=9.42cm]{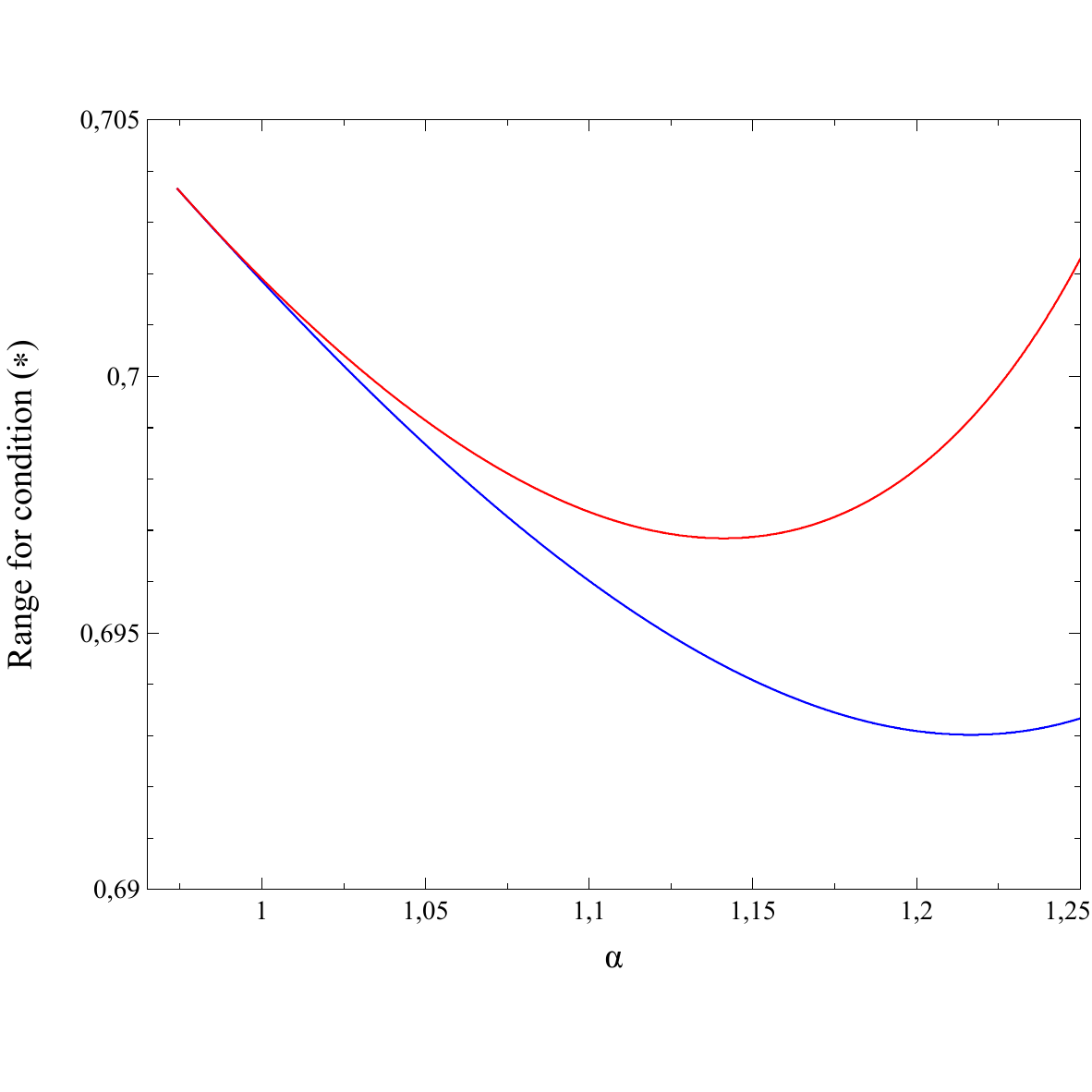}
    \caption{\label{fig:alpha} Values of $x_+(\alpha)$ (red) and $x_-(\alpha)$ (blue) between which the condition \eqref{eq:otherconditionmu} implying the stability required by \autoref{hyp:A} is satisfied. The red inset of the first plot is enlarged in the second plot.}
\end{figure}

The same procedure can be repeated for any given value of $\alpha$, so we now present numerical evidence that there exist admissible configurations for \autoref{hyp:A} for $\alpha$ in a certain range of values.
Although we can't claim this is a complete mathematical proof, the interested reader can rigorously show the analog of \autoref{thm:main} for every value of $\alpha \in (1,2)$ of their choice, by using the parameters suggested by \autoref{fig:alpha} and verifying \autoref{hyp:A}.

Given $\alpha$ we can numerically solve (by Newton's method) the second equation in \eqref{eq:intermsofxy} to find $y=y(x)$, restricting to those $x$ that satisfy the first and third items of \eqref{eq:intermsofxy}. In turn, by \eqref{eq:paramreduce} and \eqref{eq:paramreduce2}, this characterizes vorticity and position $(\xi_3,a_3)$ of the third vortex. The crucial check is now \eqref{eq:conditionmu}, which translates to the inequalities:
\begin{equation}\label{eq:otherconditionmu}\tag{*}
    c_1^2-4c_2>0,
    \quad
    2b^2-c_1-\sqrt{c_1^2-4c_2}>0.
\end{equation}
The latter are satisfied for $x$ belonging to a certain interval $(x_-(\alpha),x_+(\alpha))$, portrayed in 
\autoref{fig:alpha}.
The values of $\alpha$ at which these extremal values coincide are approximately $\alpha_-\simeq 0.97424 $ and $\alpha_+\simeq 2.13903$. For $\alpha$ bewteen those values, there exist (infinitely many) values of $x$ satisfying the above conditions, namely \autoref{hyp:A} holds true. 
\autoref{fig:alpha} reports results for $\alpha$ varying between $\alpha_\pm$ with mesh $10^{-4}$ and $x\in (0,1)$ with mesh $10^{-6}$.





\begin{acknowledgements}
The authors thank Silvia Morlacchi for the \emph{Veusz} data plot in \autoref{fig:alpha}.
    F. G. acknowledges support from the MUR Excellence Department Project awarded to the Department of Mathematics, University of Pisa, CUP I57G22000700001. 
    U.P. has received funding from the European Research Council (ERC) under the European Union’s Horizon 2020 research and innovation programme (grant agreement No. 949981).    
    U.P acknowledges the kind hospitality of Scuola Normale Superiore, where part of this work has been done.
\end{acknowledgements}

\bibliography{biblio}{}

\begin{thebibliography}{10}

\bibitem{Badin2018}
Gualtiero Badin and Anna~M. Barry.
\newblock Collapse of generalized euler and surface quasigeostrophic point
  vortices.
\newblock {\em Phys. Rev. E}, 98:023110, Aug 2018.

\bibitem{Bo99}
Thierry Bodineau and Alice Guionnet.
\newblock About the stationary states of vortex systems.
\newblock {\em Ann. Inst. H. Poincar\'e{} Probab. Statist.}, 35(2):205--237,
  1999.

\bibitem{Ca92}
E.~Caglioti, P.-L. Lions, C.~Marchioro, and M.~Pulvirenti.
\newblock A special class of stationary flows for two-dimensional {E}uler
  equations: a statistical mechanics description.
\newblock {\em Comm. Math. Phys.}, 143(3):501--525, 1992.

\bibitem{chen2024sufficient}
Jiahe Chen and Qihuai Liu.
\newblock Sufficient and necessary conditions for self-similar motions of three
  point vortices in generalized fluid systems.
\newblock {\em Physica D: Nonlinear Phenomena}, 470:134392, 2024.

\bibitem{donati2024existence}
Dimitri Cobb, Martin Donati, and Ludovic Godard-Cadillac.
\newblock Existence and uniqueness for the sqg vortex-wave system when the
  vorticity is constant near the point-vortex.
\newblock {\em arXiv preprint arXiv:2401.02728}, 2024.

\bibitem{donati2023holder}
Martin Donati and Ludovic Godard-Cadillac.
\newblock H{\"o}lder regularity for collapses of point-vortices.
\newblock {\em Nonlinearity}, 36(11):5773, 2023.

\bibitem{godard2021vortex}
Ludovic Godard-Cadillac.
\newblock Vortex collapses for the euler and quasi-geostrophic models.
\newblock {\em arXiv preprint arXiv:2101.11258}, 2021.

\bibitem{godard2023holder}
Ludovic Godard-Cadillac.
\newblock H{\"o}lder estimate for the 3 point-vortex problem with alpha-models.
\newblock {\em Comptes Rendus. Math{\'e}matique}, 361(G1):355--362, 2023.

\bibitem{GrLuRo24}
Francesco Grotto, Eliseo Luongo, and Marco Romito.
\newblock Gibbs equilibrium fluctuations of point vortex dynamics.
\newblock {\em Ann. Appl. Probab.}, 34(6):5426--5461, 2024.

\bibitem{ARMAa}
Francesco Grotto and Umberto Pappalettera.
\newblock Burst of point vortices and non-uniqueness of 2{D} {E}uler equations.
\newblock {\em Arch. Ration. Mech. Anal.}, 245(1):89--126, 2022.

\bibitem{GrRo20}
Francesco Grotto and Marco Romito.
\newblock A central limit theorem for {G}ibbsian invariant measures of 2{D}
  {E}uler equations.
\newblock {\em Comm. Math. Phys.}, 376(3):2197--2228, 2020.

\bibitem{Helmholtz1858}
H.~Helmholtz.
\newblock \"{U}ber {I}ntegrale der hydrodynamischen {G}leichungen, welche den
  {W}irbelbewegungen entsprechen.
\newblock {\em Journal f\"{u}r die Reine und Angewandte Mathematik. [Crelle's
  Journal]}, 55:25--55, 1858.

\bibitem{IwWa10}
Takahiro Iwayama and Takeshi Watanabe.
\newblock Green's function for a generalized two-dimensional fluid.
\newblock {\em Phys. Rev. E}, 82:036307, Sep 2010.

\bibitem{Ki12}
Michael K.-H. Kiessling and Yu~Wang.
\newblock Onsager's ensemble for point vortices with random circulations on the
  sphere.
\newblock {\em J. Stat. Phys.}, 148(5):896--932, 2012.

\bibitem{MaPu83}
C.~Marchioro and M.~Pulvirenti.
\newblock Euler evolution for singular initial data and vortex theory.
\newblock {\em Comm. Math. Phys.}, 91(4):563--572, 1983.

\bibitem{Marchioro88}
Carlo Marchioro.
\newblock Euler evolution for singular initial data and vortex theory: a global
  solution.
\newblock {\em Comm. Math. Phys.}, 116(1):45--55, 1988.

\bibitem{MaPu94}
Carlo Marchioro and Mario Pulvirenti.
\newblock {\em Mathematical theory of incompressible nonviscous fluids},
  volume~96 of {\em Applied Mathematical Sciences}.
\newblock Springer-Verlag, New York, 1994.

\bibitem{Reinaud21}
Jean~N. Reinaud.
\newblock Self-similar collapse of three geophysical vortices.
\newblock {\em Geophys. Astrophys. Fluid Dyn.}, 115(4):369--392, 2021.

\bibitem{vishik2018instability}
Misha Vishik.
\newblock {Instability and non-uniqueness in the Cauchy problem for the Euler
  equations of an ideal incompressible fluid. Part I and Part II}.
\newblock {\em arXiv:1805.09426, arXiv:1805.09440}, 2018.

\bibitem{zbarsky}
Samuel Zbarsky.
\newblock From point vortices to vortex patches in self-similar expanding
  configurations.
\newblock {\em Comm. Math. Phys.}, 388(2):707--733, 2021.

\end{thebibliography}
\bibliographystyle{plain}

\end{document}